\documentclass[12pt]{article}
\usepackage{diagbox}
\usepackage{mathtools}
\usepackage{bbm}
\usepackage{amsmath}
\usepackage{hyperref}
\allowdisplaybreaks

\mathtoolsset{showonlyrefs}

\renewcommand{\emph}[1]{\textit{#1}}
\usepackage{enumerate,amsmath,amsthm,latexsym,amssymb}
\usepackage{color}\usepackage{graphicx}

\def\tu{\tilde{u}}

\definecolor{brown}{cmyk}{0, 0.72, 1, 0.45}
\definecolor{grey}{gray}{0.5}

\newcommand{\old}[1]{}

\newcounter{rot}%\addtocounter{rot}{1}, \therot

\newcommand{\ignore}[1]{}

\newcommand{\set}[1]{\left\{#1\right\}}

\def\cP{\mathcal{P}}

\newcommand{\proofend}{\hspace*{\fill}\mbox{$\Box$}\\ \medskip\\ \medskip}

\def\ii_(#1,#2){i_{#1}^{#2}}

\def\a{\alpha}

\def\D{\Delta}
\def\e{\varepsilon}

\def\G{\Gamma}

\def\P{\Pi}

\def\s{\sigma}
\def\t{\tau}

\def\1{{\bf 1}}
\def\0{{\bf 0}}

\newcommand{\brac}[1]{\left( #1 \right)}

\renewcommand{\Pr}{\operatorname{\bf Pr}}

\newtheorem{theorem}{Theorem}[section]

\newtheorem{lemma}[theorem]{Lemma}

\newtheorem{remthm}[theorem]{Remark}

\newcounter{thmtemp}

\newcommand{\nospace}[1]{}

\def\path{\operatorname{PATH}}

\newcommand{\beq}[2]{\begin{equation}\label{#1}#2\end{equation}}

\def\hd{\widehat{d}}
\def\hp{\widehat{p}}
\def\hG{\widehat{\G}}

\author{Michael Anastos, Alan Frieze\thanks{Research supported in part by NSF grant DMS1362785
} and Wesley Pegden\thanks{Research supported in part by NSF grant DMS1363136}\\Department of Mathematical Sciences\\Carnegie Mellon University\\Pittsburgh PA 15213} 
\begin{document}
\title{Constraining the clustering transition for\\ colorings of sparse random graphs}
\maketitle
\begin{abstract}
Let $\Omega_q$ denote the set of proper $q$-colorings of the random graph $G_{n,m}, m=dn/2$ and let $H_q$ be the graph with vertex set $\Omega_q$ and an edge $\set{\s,\t}$ where $\s,\t$ are mappings $[n]\to[q]$ iff $h(\s,\t)=1$. Here $h(\s,\t)$ is the Hamming distance $|\set{v\in [n]:\s(v)\neq\t(v)}|$. We show that w.h.p. $H_q$ contains a single giant component containing almost all colorings in $\Omega_q$ if $d$ is sufficiently large and $q\geq \frac{cd}{\log d}$ for a constant $c>3/2$.
\end{abstract}
\section{Introduction}\label{intro}
In this short note, we will discuss a structural property of the set $\Omega_q$ of proper $q$-colorings of the random graph $G_{n,m}$, where $m=dn/2$ for some large constant $d$.  For the sake of precision, let us define $H_q$ to be the graph with vertex set $\Omega_q$ and an edge $\set{\s,\t}$ iff $h(\s,\t)=1$ where $h(\s,\t)$ is the Hamming distance $|\set{v\in [n]:\s(v)\neq\t(v)}|$. In the Statistical Physics literature the definition of $H_q$ may be that colorings $\s,\t$ are connected by an edge in $H_q$ whenever $h(\s,\t)=o(n)$. Our theorem holds a fortiori if this is the case.

Heuristic evidence in the statistical physics literature (see for example \cite{ZK}) suggests there is a \emph{clustering transition} $c_d$ such that for $q>c_d$, the graph $H_q$ is dominated by a single connected component, while for $q<c_d$, an exponential number of components are required to cover any constant fraction of it; it may be that $c_d\approx\frac{d}{\log d}$. (Here $A(d)\approx B(d)$ is taken to mean that $A(d)/B(d)\to 1$ as $d\to \infty$. We do not assume $d\to\infty$, only that $d$ is a sufficiently large constant, independent of $n$.)  Recall that $G_{n,m}$ for $m=dn/2$ becomes $q$-colorable around $q\approx\frac{d}{2\log d}$ or equivalently when $d\approx 2q\log q$, \cite{AN,CV}.
In this note, we prove the following:
\begin{theorem}\label{th1}
If $q\geq \frac{cd}{\log d}$ for constant $c>3/2$, and $d$ is sufficiently large, then w.h.p.~$H_q$ contains a giant component that contains almost all of $\Omega_q$.
\end{theorem}
 In particular, this implies that the clustering transition $c_d$, if it exists, must satisfy $c_d\leq \frac 3 2 \frac{d}{\log d}$.  

 Theorem \ref{th1} falls into the area of ``Structural Properties of Solutions to Random Constraint Satisfaction Problems''. This is a growing area with connections to Computer Science and Theoretical Physics.  In particular, much of the research on the graph $H_q$ has been focussed on the structure near the {\em colorability threshold}, e.g. Bapst, Coja-Oghlan, Hetterich, Rassman and Vilenchik \cite{BCHRV}, or the {\em clustering threshold}, e.g. Achlioptas, Coja-Oghlan and Ricci-Tersenghi \cite{ACR}, Molloy \cite{M}.  Other papers heuristically identify a sequence of phase transitions in the structure of $H_q$, e.g., Krz\c{a}kala, Montanari, Ricci-Tersenghi, Semerijan and Zdeborov\'a \cite{KMRSZ}, Zdeborov\'a and Krz\c{a}kala \cite{ZK}.  The existence of these transitions has been shown rigorously for some other CSPs. One of the most spectacular examples is due to Ding, Sly and Sun \cite{sly} who rigorously showed the existence of a sharp satisfiability threshold for random $k$-SAT.

An obvious target for future work is improving the constant in Theorem \ref{th1} to 1. We should note that Molloy \cite{M} has shown that w.h.p. there is no giant component if $q\leq \frac{(1-\e_d)d}{\log d}$, for some $\e_d>0$. Looking in another direction, it is shown in \cite{DFFV} that w.h.p. $H_q,q\geq d+2$ is connected. This implies that Glauber Dynamics on $\Omega_q$ is ergodic. It would be of interest to know if this is true for some $q\ll d$.

Before we begin our analysis, we briefly explain the constant 3/2. We start with an arbitrary $q$-cloring and then re-color it using only approximately $\approx d/\log d$ of the given colors. We then use a disjoint set of approximately $d/2\log d$ colors to re-color it with a target $\chi\approx \frac{d}{2\log d}$ coloring $\t$.  
\section{Greedily Re-coloring}
Our main tool is a theorem from Bapst, Coja-Oghlan and Efthymiou \cite{BCE} on planted colorings. We consider two ways of generating a random coloring of a random graph. We will let $Z_q=|\Omega_q|$. The first method is to generate a random graph and then a random coloring. In the second method, we generate a random (planted) coloring and then generate a random graph compatible with this coloring.

{\bf Random coloring of the random graph $G_{n,m}$:} Here we will assume that $m$ is such that w.h.p. $Z_q>0$.
\begin{enumerate}[(a)]
\item Generate $G_{n,m}$ subject to $Z_q>0$.
\item Choose a $q$-coloring $\s$ uniformly at random from $\Omega_q$.
\item Output $\Pi_1=(G_{n,m},\s)$.
\end{enumerate}
{\bf Planted model:} 
\begin{enumerate}
\item Choose a random partition of $[n]$ into $q$ color classes $V_1,V_2,\ldots,V_q$ subject to 
\beq{sizes}{
\sum_{i=1}^q\binom{|V_i|}{2}\leq \binom{n}{2}-m.
}
\item Let $\G_{\s,m}$ be obtained by adding $m$ random edges, each with endpoints in different color classes.
\item Output $\Pi_2=(\G_{\s,m},\s)$.
\end{enumerate}
We will use the following result from \cite{BCE}:
\begin{theorem}\label{plant}
Let $d=2m/n$ and suppose that $d\leq 2(q-1)\log(q-1)$. Then $\Pr(\P_2\in \cP)=o(1)$ implies that $\Pr(\P_1\in\cP)=o(1)$ for any graph+coloring property $\cP$.
\end{theorem}
Consequently, we will use the planted model in our subsequent analysis. Let
$$q_0=\frac{q}{q-1}\cdot\frac{d}{\log d-7\log\log d}\approx \frac{d}{\log d}.$$
The property $\cP$ in question will be: ``the given $q$-coloring can be reduced via single vertex color changes to a $q_0$ coloring'' where $\a>1$ is constant.

In a random partition of $[n]$ into $q$ parts, the size of each part is distributed as $Bin(n,q^{-1})$ and so the Chernoff bounds imply that w.h.p. in a random partition each part has size $\frac{n}{q}\brac{1\pm \frac{\log n}{n^{1/2}}}$.

We let $\G$ be obtained by taking a random partition $V_1,V_2,\ldots,V_q$ and then adding $m=\frac12dn$ random edges so that each part is an independent set. These edges will be chosen from $$N_q=\binom{n}{2}-(1+o(1))q\binom{n/q}{2}=(1-o(1))\frac{n^2}{2}\brac{1-\frac1q}$$ 
possibilities. So, let $\hd=\frac{mn}{N_q}\approx\frac{dq}{q-1}$ and replace $\G$ by $\hG$ where each edge not contained in a $V_i$ is included independently with probability $\hp=\frac{\hd}{n}$. $V_1,V_2,\ldots,V_q$ constitutes a coloring which we will denote by $\s$. Now $\hG$ has $m$ edges with probability $\Omega(n^{-1/2})$ and one can check that the properties required in Lemmas \ref{indsparse0} and \ref{indsparse} below all occur with probability $1-o(n^{-1/2})$ and so we can equally well work with $\hG$.

Now consider the following algorithm for going from $\s$ via a path in $\Omega_q$ to a coloring with significantly fewer colors. It is basically the standard greedy coloring algorithm, as seen in Bollob\'as and Erd\H{o}s \cite{BE}, Grimmett and McDiarmid \cite{GM} and in particular Shamir and Upfal \cite{SU} for sparse graphs. 

In words, it goes as follows. At each stage of the algorithm, $U$ denotes the set of vertices that have not been re-colored. Having used $r-1$ colors to color some subset of vertices we start using color $r$. We let $W_j=V_j\cap U$ denote the uncolored vertices of $V_j$ for $j\geq 1$. We then let $k$ be the smallest index $j$ for which $W_j\neq\emptyset$. This is an independent set and so we can re-color the vertices of $W_k$, one by one, with the color $r$. We let $U_r\subseteq U$ denote the set of vertices that may possibly be re-colored $r$ by the algorithm i.e. those vertices with no neighbors in $C_r$, the current set of vertices colored $r$.  Each time we re-color a vertex with color $r$, we remove its neighbors from $U_r$. We continue with color $r$, until $U_r=\emptyset$. After which, $C_r$ will be the set of vertices that are finally colored with color $r$.

At any stage of the algorithm, $U$ is the set of vertices whose colors have not been altered. The value of $L$ in line D is $n/\log^2\hd$. 
\begin{tabbing}
{\sc algorithm greedy re-color}\\
{\bf beg}\={\bf in}\\
\>Initialise: $r=0,U=[n], C_0\gets \emptyset$;\=\\
\>{\bf repeat};\\
\>$r\gets$\=$r+1$, $C_r\gets\emptyset$;\\
\>\>Let $W_j=V_j\cap U$ for $j\geq 1$ and let $k=\min\set{j:W_j\neq \emptyset}$;\\
{\bf A:}\>\>$C_r\gets W_k,U\gets U\setminus C_r,U_r\gets U\setminus\set{\text{neighbors of $C_r$ in $\hG$}}$;\\
\>\>If $r<k$, re-color every vertex in $C_r$ with color $r$;\\
{\bf B:}\>\>{\bf rep}\={\bf eat} \=(Re-color some more vertices with color $r$);\\
{\bf C:}\>\>\>Arbitrarily choose $v\in U_r$, $C_r\gets C_r+v$, $U_r\gets U_r-v$;\\
\>\>\>$U_r\gets U_r\setminus\set{\text{neighbors of $v$ in $\hG$}}$;\\ 
\>\>{\bf until} $U_r=\emptyset$;\\
{\bf D:}\>{\bf until} $|U|\leq L$;\\
\>Re-color $U$ with $\frac{\hd}{\log^2\hd}+2$ unused colors from our initial set of $q_0$ colors;\\
{\bf end}
\end{tabbing}
We first observe that each re-coloring of a singe vertex $v$ vertex in line C can be interpreted as moving from a coloring of $\Omega_q$ to a neighboring coloring in $H_q$. This requires us to argue that the re-coloring by {\sc greedy re-color} is such that the coloring of $\hG$ is proper at all times. We argue by induction on $r$ that the coloring at line A is proper. When $r=1$ there have been no re-colorings. Also, during the loop beginning at line B we only re-color vertices with color $r$ if they are not neighbors of the set $U_r$ of vertices colored $r$. This guarantees that the coloring remains proper until we reach line D. The following lemma shows that we can then reason as in Lemma 2 of Dyer, Flaxman, Frieze and Vigoda \cite{DFFV}, as will be explained subsequently. 
\begin{lemma}\label{indsparse0}
Let $p=m/\binom{n}{2}=\D/n$ where $\D$ is some sufficiently large constant.
With probability $1-o(n^{-1/2})$, every $S\subseteq [n]$ with $s=|S|\leq n/\log^2\D$ contains at most $s\D/\log^2\D$ edges. 
\end{lemma}
The above lemma, is Lemma 7.7(i) of  Janson, {\L}uczak and Ruci\'nski \cite{JLR} and it implies that if $\D=\hd$ then w.h.p. $\hG_U$  at line D contains no $K$-core, $K=\frac{2\hd}{\log^2\hd}+1$. Here $\hG_U$ denotes the sub-graph of $\hG$ induced by the vertices $U$. For a graph $G=(V,E)$ and $K\geq 0$, the $K$-core is the unique maximal set $S\subseteq V$ such that the induced subgraph on $S$ has minimum degree at least $K$. A graph without a $K$-core is {\em $K$-degenerate} i.e. its vertices can be ordered as $v_1,v_2,\ldots,v_n$ so that $v_i$ has at most $K-1$ neighbors in $\{v_1,v_2,\ldots,v_{i-1}\}$. To see this, let $v_n$ be a vertex of minimum degree and then apply induction.
 
We argue now that we can re-color the vertices in $U$ with $K+1$ new colors, all the time following some path in $H_q$. Let $v_1,\dots,v_n$ denote an ordering of $U$ such that
the degree of $v_i$ is less than $K$ in the subgraph $\hG_i$ of $\hG$ induced by $\{v_1,v_2,\dots,v_i\}$. We will prove the claim by induction. The claim is trivial for $i=1$. By induction there is a path $\s_0,\s_1,\ldots,\s_r$ from the coloring $\s_0$ of $U$ at line B, restricted to $\hG_{i-1}$ using only $K+1$ colors to do the re-coloring.

Let $(w_j,c_j)$ denote the $(vertex,color)$ change defining the edge $\set{\s_{j-1},\s_j}$.  We construct a path (of length $\leq 2r$) that re-colors $\hG_i$. For $j=1,2,\ldots,r$, we will re-color $w_j$ to color $c_j$, if no neighbor of $w_j$ has color $c_j$. Failing this, $v_i$ must be the only neighbor of $w_j$ that is colored $c_j$. This is because $\s_r$ is a proper coloring of $\hG_{i-1}$. Since $v_i$ has degree less than $K$ in $\hG_i$, there exists a new color for $v_i$ which does not appear in its neighborhood. Thus, we first re-color $v_i$ to any new (valid) color, and then we re-color $w_j$ to $c_j$, completing the inductive step. Note that because the colors used in Step D have not been used in Steps A,B,C, this re-coloring does not conflict with any of the coloring done in Steps A,B,C.

We need to show next that each Loop B re-colors a large number of vertices. Let $\a_1(G)$ denote the minimim size of a {\em maximal} independent set of a graph $G$ i.e. an independent set that is not contained in any larger independent set. The round will re-color at least $\a_1(\G_U)$ vertices, where $U$ is as at the start of Loop B. The following result is from Lemma 7.8(i) of \cite{JLR}.
\begin{lemma}\label{indsparse}
Let $p=m/\binom{n}{2}=\D/n$ where $\D$ is some sufficiently large constant.
$\a_1(G_{n,m})\geq \frac{\log \D-3\log\log \D}{p}$ with probability $1-o(n^{-1/2})$. (see Lemma 7.8(i)).
\end{lemma}
Suppose now that we take $u_0$ to be the size of $U$ at the beginning of Step A and that $u_t$ is the size of $U$ after $t$ vertices have been finally colored $r$. Thus we assume that $u_{|W_k|}$ is the size of $U$ at the start of Step B. We observe that, 
\beq{uchange}{
u_{t+1}\text{ stochastically dominates } u_t-Bin(u_t,\hp)-1.
}
This is because the edges inside $U$ are unconditioned by the algorithm and because $v\in V_j$ has no neighbors in $V_j$ for $j\geq 1$. On the other hand, if we apply Algorithm {\sc greedy re-color} to $G_{n,\hp}$ then \eqref{uchange} is replaced by the recurrence
\beq{uchange1}{
\tu_{t+1}= \tu_t-Bin(\tu_t,\hp)-1.
}
(Putting $V_j=\set{j}$ means that {\sc greedy re-color} is running on $G_{n,\hp}$.)

Comparing \eqref{uchange} and \eqref{uchange1} we see that we can couple the two applications of {\sc greedy re-color} so that $u_t\geq \tu_t$ for $t\geq 0$. Now the application of Loop B re-colors a maximal independent set of the graph $\hG_U$ induced by $U$ as it stands at the beginning of the loop. The size of this set dominates the size of a maximal independent set in the random graph $G_{|U|,p}$.  So if we generate $G_{|U|,p}$ and then delete some edges, we see that every independent set of $G_{|U|,p}$ will be contained in an independent set of $\G_U$. And so using Lemma \ref{indsparse} we see that w.h.p. each execution of Loop B re-colors at least 
$$\frac{\log(\hd/\log^2\hd)-3\log\log(\hd/\log^2\hd)}{\hd}n\geq \frac{q-1}{q} \cdot\frac{\log d-6\log\log d}{d}n$$ 
vertices, for $d$ sufficiently large. We have replaced $\D$ of Lemma \ref{indsparse} by $\hd/\log^2\hd$ to allow for the fact that we hae replaced $n$ by $|U|\geq L$. Consequently, at the end of Algorithm {\sc greedy re-color} we will have used at most 

\beq{}{
\frac{q}{q-1}\cdot\frac{d}{\log d-6\log\log d}+\frac{\hd}{\log^2\hd}+2 \leq \frac{q}{q-1}\cdot\frac{d}{\log d-7\log\log d}=q_0
}
colors. The term $\frac{\hd}{\log^2\hd}+2$ arises from the re-coloring of $U$ at line D.

{\bf Finishing the proof:} Now suppose that $q\geq \frac{cd}{\log d}$ where $d$ is large and $c>3/2$. Fix a particular $\chi$-coloring $\tau$. We prove that almost every $q$-coloring $\s$ can be transformed into $\t$ changing one color at a time. It follows that for almost every pair of $q$-colorings $\s$, $\s'$ we can transform $\s$ into $\s'$ by first transforming $\s$ to $\t$
and then reversing the path from $\s'$ to $\t$.  

We proceed as follows. The algorithm {\sc greedy re-color} takes as input: (i) the coloring $\s$ and (ii) a specific subset of $q_0$ colors from $\set{1, ..., q}$ that are not used in $\t$. W.h.p. it transforms the input coloring into a coloring using only those $q_0$ colors. Then we process the color classes of $\t$, re-coloring vertices to their $\t$-color. When we process a color class $C$ of $\t$, we switch the color of vertices in $C$ to their $\t$-color $i_C$ one vertex at a time. We can do this because when we re-color a vertex $v$, a neighbor $w$ will currently either have one of the $q_0$ colors used by {\sc greedy re-color} and these are distinct from $i_C$. Or $w$ will have already been been re-colored with its $\t$-color which will not be color $i_C$. This  proves Theorem \ref{th1}.
\proofend

\end{document}